\documentclass[13pt,oneside,english,thmsb]{article}

\usepackage[T1]{fontenc}
\usepackage[latin1]{inputenc}
\usepackage{amssymb}

\usepackage{float}
\usepackage{amsmath}

\usepackage{graphicx}
\usepackage{latexsym}
\newcommand{\nbf}{\mathbf{n}}
\newcommand{\ibf}{\mathbf{i}}
\newcommand{\jbf}{\mathbf{j}}
\newcommand{\Obf}{\text{\mathversion{bold}{$\mathcal{O}_\mathbf{n}$}}}
\newcommand{\Ibf}{\text{\mathversion{bold}{$\mathcal{I}_\mathbf{n}$}}}
\newcommand{\Pbf}{\text{\mathversion{bold}{$P$}}}
\newcommand{\Ebf}{\text{\mathversion{bold}{$E$}}}

\newcommand{\argmin}{\operatornamewithlimits{argmin}}
\newtheorem{theo}{Theorem}
\newtheorem{prop}{Proposition}
\newtheorem{lem}{Lemma}
\newtheorem{cor}{Corollary}

\def\1{\hbox{\rm 1\hskip -3pt I}}

\title{Robust quantile estimation and prediction for spatial processes}
 \author{Sophie Dabo-Niang \thanks{Corresponding author: sophie.dabo@univ-lille3.fr}\and Baba Thiam
\thanks{ Laboratoire EQUIPPE, Université Charles-De-Gaulle, Lille
3, Maison de la Recherche, domaine universitaire du Pont de Bois, BP 60149, 59653 Villeneuve d'Ascq cedex, France. baba.thiam@univ-lille3.fr}}

\begin{document}

\maketitle
\noindent

\begin{abstract}
In this paper, we present a
statistical framework for modeling conditional quantiles of spatial processes assumed to
be strongly mixing in space. We establish the $L_1$ consistency
and the asymptotic normality of the kernel conditional quantile estimator in the case
of random fields. We also define a nonparametric spatial predictor and illustrate the methodology used with some simulations. \\

\end{abstract}

\section{Introduction}\label{section1}
Let $(X,Y)$ be a pair of random variables with values in $\mathbb{R}^d\times \mathbb{R}$ and defined on a
probability space $(\Omega,\mathcal{A},\boldmath{P})$. Assume that the joint density of $(X,Y)$ and the marginal density of $X$ exist and are denoted respectively by $f(x,y)$ and $g(x)$. In the following, we suppose that $F(\cdot|x)$, the conditional distribution function of $Y$ given $X=x$ exists and we denote by $f(\cdot|x)$ the density of $Y$ given $X=x$.  For $p \in ]0,1[$ and for fixed $x\in \mathbb{R}^d$, let $\mu_p(x)$ be the conditional quantile of order $p$ of $F(\cdot|x)$, that can be seen as a solution of the equation $F(y|x)=p$.
Another alternative characterization of the $p$$th$ conditional quantile (see for example Gannoun et al.
\cite{gannoun}) is $\mu_p(x)= \argmin_{\theta \in \mathbb{R}}\Ebf\left[(2p-1)(Y-\theta)+|Y-\theta|\ \vert X=x\right]$.\\
We are interested to the non-parametric estimation of $\mu_p(x)$ in the case of spatial dependent observations.
Nonparametric conditional quantile estimation technics have already been developed for non spatial (independent or mixing) real valued processes. Such results have
provided useful tools for solving for example some prediction problems of strictly stationary processes satisfying the $\alpha$-mixing condition. The existing results in the non-spatial case include the works of Matzner-L\o ber \cite{matzner}, Collomb \cite{collomb}, Gannoun et al. \cite{gannoun}, Laksaci et al. \cite{laksaci}. \\
In nonparametric spatial estimation, the existing works concern mainly the estimation of a probability density and regression functions, see the key references: Tran \cite{tran}, Biau and Cadre \cite{biau}, Carbon et al. \cite{carbon}.\\
For the spatial quantile conditional estimation case, there exist only few results in our knowledge. Abdi et al. \cite{abdi} considered the pointwise $p-$mean and almost complete consistencies of a double kernel quantile estimator for real-valued random fields.  Hallin {\it et al.} \cite{hallin} give a Bahadur representation and asymptotic normality results of the local linear quantile estimator. Laksaci and Fouzia \cite{laksaci-fouzia} consider the case where the regressor take their values in a semi-metric space and show the strong and weak consistency of the conditional quantile. \\
In this paper, we will go beyond all these last spatial works and provide the $L_1$ consistency and an asymptotic normality  of a kernel conditional quantile estimate of a strictly stationary spatial process satisfying the
$\alpha$-mixing condition. In addition, we employe our results to
solve some nonparametric prediction problems. The organization of this paper is as follows. The estimation procedure is presented in Section 2. Section 3 gives some necessary conditions and then establishes the
main asymptotic results. Section 4 is devoted to simulations results. Technical proofs are given in Section 5.
\section{Nonparametric estimator of the conditional quantile}\label{section_cons}
Let us consider a strictly stationary process $\left((X_\ibf,Y_\ibf), \ibf \in \Ibf\right)$ with values in $\mathbb{R}^d\times \mathbb{R}$ where $(X_\ibf,Y_\ibf)$ has the same distribution as $(X,Y)$.
For $\nbf\in (\mathbb{N^*})^N$, we define a rectangular region $\Ibf$ by $\Ibf=\{\ibf=(i_1,\ldots,i_N) \in (\mathbb{N^*})^N, \ \  1\leq i_k\leq n_k, \ \
k=1,\ldots,N\}$. We set $\widehat{\nbf}=n_1\ldots n_N$, and we write $\nbf \to \infty$ if $\min_{k=1,\ldots,N}n_k\to \infty$. The well known kernel estimates of $f$ and $g$
are defined by

$\displaystyle f_\nbf(x,y)  =  \frac{1}{\widehat\nbf h^{d+1}}\sum_{\ibf\in\Ibf}K\left(\frac{x-X_\ibf}{h}\right)w\left(\frac{y-Y_\ibf}{h}\right)$ and $\displaystyle g_\nbf(x)  = \frac{1}{\widehat\nbf h^d}\sum_{\ibf\in\Ibf}K\left(\frac{x-X_\ibf}{h}\right)$,
where $K$ and $w$ are two probability density functions, and the bandwidths $h=h(\nbf)$ is a sequence of positive real numbers such that $h \to 0$ as $\nbf \to \infty$.
The kernel estimate of the conditional density $f_\nbf(y|x)$ is naturally defined by the ratio $f_\nbf(x,y)$ over $g_\nbf(x)$ while the estimator of the conditional distribution function (see the one introduced by Roussas \cite{roussas}) is defined by
\begin{eqnarray*}
F_{\nbf}(y|x)  &= & \frac{\psi_{\nbf}(x,y)}{g_\nbf(x)}\1_{\left\{g_\nbf(x)\neq 0\right\}},\\
\psi_{\nbf}(x,y) & = & \frac{1}{\widehat\nbf h^{d+1}}\sum_{\ibf\in\Ibf}K\left(\frac{x-X_\ibf}{h}\right)\int_{-\infty}^yw\left(\frac{z-Y_\ibf}{h}\right)dz.
\end{eqnarray*}
For a fixed $x$, the estimator of the $pth$ conditional quantile noted $\mu_{p,\nbf}(x)$  can be defined as the root of the equation $ F_{\nbf}(z|x)=p.$
 Alternatively, one can consider the local constant estimator defined by
\begin{eqnarray*}\label{local_estimator}
\nu_{p,\nbf}(x)=\argmin_{\theta\in\mathbb{R}}\sum_{\ibf\in\Ibf}\left(|Y_\ibf-\theta|+(2p-1)(Y_\ibf-\theta)\right)K\left(\frac{x-X_\ibf}{h}\right).
\end{eqnarray*}
In this paper, we will focus on the study of the asymptotic behavior of $\mu_{p,\nbf}$. For the study of $\nu_{p,\nbf}$, one can adapt the technics developed in Zhou \cite{zhou}.

\section{Main results}\label{results}
To establish the asymptotic results, we will suppose that the sequence $\left(X_\ibf,Y_\ibf)\right)_{\ibf \in (\mathbb{N}^*)^N}$ satisfies the following
mixing condition: there exists a function $\chi:\mathbb{R}^+\to \mathbb{R}^+$ with $\chi(t)\downarrow 0$ as $t \to \infty$, such that whenever $E,E'\subset
(\mathbb{N}^*)^N$ with finite cardinals,
\begin{eqnarray*}
\alpha\left(\mathcal{B}(E),\mathcal{B}(E')\right) & :=  & \sup\left\{|\Pbf(A\cap B)-\Pbf(A)\Pbf(B)|; \ A \in \mathcal{B}(E), B \in \mathcal{B}(E') \right\}\\
& \leq & \phi(Card E,Card E')\chi(dist(E,E')),
\end{eqnarray*}
where $\mathcal{B}(E)$ (resp. $\mathcal{B}(E')$) denotes the Borel $\sigma$-fields generated by $(X_\ibf,Y_\ibf)_{\ibf \in E}$ (resp. $(X_\ibf,Y_\ibf)_{\ibf \in E'}$), Card $E$ (resp. Card $E'$) the cardinality of $E$ (resp. $E'$), $dist(E,E')$ the Euclidean distance between $E$ and $E'$, and $\phi:\mathbb{N}^2\to \mathbb{R}^+$ is a symmetric positive function which is non decreasing in each variable. Throughout this paper, we will assume that $\phi$ satisfies
\begin{equation}
\phi\left(n,m\right)\leq C\min\left(n,m\right),\,\,\,\,\,\,\forall
n,m\in\mathbb{N}\label{mixing1}\end{equation}
 or \begin{equation}
\phi\left(n,m\right)\leq
C\left(n+m+1\right)^{\kappa},\,\,\,\,\,\,\,\,\,\forall
n,m\in\mathbb{N}\label{mixing2}\end{equation}
 for some $\kappa\geq1$ and some $C>0$.
If $\phi\equiv 1$, then the field $(X_\ibf,Y_\ibf)_{\ibf\in
(\mathbb{N}^*)^{\mathbb{N}}}$ is called \textit{strongly mixing}. In
this paper, we consider the case where $\chi (i)$ tends to zero at a
polynomial rate, that is,
\begin{eqnarray}
\label{mel}
\chi(i)= O(i^{-\beta}),
\end{eqnarray}
with $\beta>0$. %As mentioned in Carbon et al. (1997), if \eqref{mel} holds for $\beta>2N$, then
%\begin{eqnarray}
%\label{serie}
%\sum_{i=1}^\infty i^{N-1}\left(\chi(i)\right)^a<\infty
%\end{eqnarray}
%for some $0<a<1/2$.
We fix a compact subset $S$ of $\mathbb{R}^d$. Denote $a=\inf\{y: F(y|x)>0 \}$ and $b=\sup\{y: F(y|x)<1\}$, we will suppose that $\mathcal{V}\subseteq[a,b]$ is a compact neighborhood of the unknown quantile $\mu(x)$. For mixing coefficients with polynomial decreasing rate \eqref{mel}, the constraints on the bandwidth will be related to $\beta$ by means of
 \[
\theta_{1}=\frac{N(d+1)(d+2)+(d+1)\beta}{\beta-N(d+5)},\ \ \theta_{2}=\frac{N(d+2)-\beta}{\beta-N(d+5)},\]

\[
\theta_{3}=\frac{N(d^2+4d+2)+(d+1)\beta}{\beta-N(d+4+2\kappa)},\ \ \theta_{4}=\frac{N(d+1)-\beta}{\beta-N(d+4+2\kappa)}.\]
Denote $\Omega_\nbf=\sqrt{\frac{\log \widehat\nbf}{\widehat\nbf h^d}}$.
Let $\varepsilon$ be an arbitrary small positive number and set
$u(\nbf)=\prod_{i=1}^N(\log n_i)(\log\log n_i)^{1+\varepsilon}$. It
is clear that $\sum_{\nbf \in \mathbb{Z}^N}1/(\widehat\nbf
u(\nbf))<\infty$.
In the sequel, we use the following hypotheses.
\begin{description}
\item(A1) $f$ and $g$ are respectively continuous on $\mathbb{R}^{d+1}$ and $\mathbb{R}^d$, $g$ satisfies a Lipschitz condition, $g(x)>0, \forall x\in S$.
\item(A2) There exists $D\geq 0$ such that the pairs $(X_\ibf,X_\jbf)$ and $\left((X_\ibf,Y_\ibf), (X_\jbf,Y_\jbf)\right)$ admit a density, say $g_{\ibf,\jbf}$ and $f_{\ibf,\jbf}$, as soon as $dist(\ibf,\jbf)>D$. Moreover, for some constant $c\geq 0$,
\begin{eqnarray*}
|f_{\ibf,\jbf}(s,t)-f(s)f(t)|\leq c, \ \ \forall s,t\in \mathbb{R}^{d+1}
\end{eqnarray*}
and
\begin{eqnarray*}
|g_{\ibf,\jbf}(u,v)-g(u)g(v)|\leq c, \ \ \forall u,v\in \mathbb{R}^{d}.
\end{eqnarray*}
\item(A3) i) $f^{(i,j)}(\cdot,\cdot)=\frac{\partial^{i+j}f}{\partial x_l^i\partial y^j}(\cdot,\cdot)$ exists, is bounded and integrable for $0\leq i+j\leq 2$, and $0\leq l\leq d$.\\
         ii) $F(y|x)$ has continuous second partial derivatives with respect to $x$.
\item (A4) $F(y|x)$ has continuous second derivative with respect to $y$.
\item(A5) The kernel $K$ is integrable, symmetric and is a lipschitzian density function on $\mathbb{R}^d$ with compact support. Moreover $\int_{\mathbb{R}^d}\|s\|^2K(s)ds<\infty$.
\item(A6) The kernel $w$ is a symmetric and lipschitzian density function on $\mathbb{R}$ and has compact support.
%and satisfies the condition $|v|w(v)\to 0$ when $v\to \infty$ and $\int_{\mathbb{R}}s^2w(s)ds<\infty$.
\item (A7) $\lim_{\nbf \to \infty}\widehat\nbf h^{d+2}\left(\log \widehat\nbf\right)^{-1}= 0$.
\item (A8) The function $\mu_p(x)$ satisfies a uniform uniqueness property on $S$:
\begin{eqnarray*}
\forall \varepsilon >0, \exists \eta>0, \forall r:S\to \mathbb{R},
\sup_{x \in S}|\mu_p(x)-r(x)|\geq \varepsilon \Rightarrow \sup_{x \in
S}\left|F(\mu_p(x)|x)-F(r(x)|x)\right|\geq \eta.
\end{eqnarray*}
\item (A9) $\widehat{\mathbf{n}}h^{\theta_{1}}(\log\widehat{\mathbf{n}})^{\theta_{2}}
\left(u(\mathbf{n})\right)^{\frac{-2N}{\beta-N(d+5)}}\rightarrow\infty$.
\item (A10) $\widehat{\mathbf{n}}h^{\theta_{3}}(\log\widehat{\mathbf{n}})^{\theta_{4}}
\left(u(\mathbf{n})\right)^{\frac{-2N}{\beta-N(d+4+2\kappa)}}\rightarrow\infty.$

\end{description}

\medskip

\noindent\textbf{Comments on the hypotheses:} \\
\noindent Assumptions (A9) and (A10) imply conditions (3.7) and (3.8) of Theorem 3.3 in Carbon et al. \cite{carbon} and they also imply the classical condition $\widehat\nbf h^{d+1}/\log\widehat\nbf \to \infty$.\\
Assumption (A8) is introduced for getting consistency results on the quantile from those of the conditional distribution.\\

\noindent In order to state the asymptotic results, we will suppose that (A9) and (\ref{mixing1}) or (A10) and (\ref{mixing2}) are satisfied. The following two theorems give uniform almost sure convergence results of respectively $F_{\nbf}(y|x)$ and $\mu_{p,\nbf}(x)$ and permit to establish the $L_1$ consistency of $\mu_{p,\nbf}(x)$ (see Corollary \ref{convl1}).
\begin{theo}\label{convcondps}
Assume (A1)-(A7) hold, then
\begin{eqnarray*}
\sup_{y \in \mathcal{V}}\sup_{x\in S}|F_{\nbf}(y|x)-F(y|x)| &= & O\left(\Omega_\nbf\right)\ \ \mbox{a.s.}
\end{eqnarray*}
\end{theo}

\begin{theo}\label{convmu}
If (A1)-(A8) are satisfied, then we have
\begin{eqnarray*}
\sup_{x\in S}|\mu_{p,\nbf}(x)-\mu_p(x)|\stackrel{a.s.}{\to} 0.
\end{eqnarray*}
\end{theo}
\begin{cor}\label{convl1}
Assume (A1)-(A8) hold, then
\begin{eqnarray*}
\Ebf\left[\{\mu_{p,\nbf}(X_\nbf)-\mu_p(X_\nbf)\}\1_{\{X_\nbf\in
S\}}\right]\stackrel{a.s.}{\to} 0.
\end{eqnarray*}
\end{cor}

\medskip

%\section{Asymptotic normality}
\noindent To establish the following asymptotic normality of $\mu_{p,\nbf}(x)$ (Theorem \ref{normquantile1}), we will suppose that for any $(x,y)\in S\times\mathcal{V}$, there exists $c>0$ such that $f(y|x)>c$. Moreover, we will assume
that the following additional conditions on the bandwidth hold for some $0<\gamma<1$.
\begin{description}
\item (C1) $\widehat{\mathbf{n}}h^{d(1+2N(1-\gamma))}\to \infty$.
\item (C2) There exists a sequence of positive integers $q=q_\nbf\to \infty$ with $q=o\left((\widehat{\mathbf{n}}h^{d(1+2N(1-\gamma))})^{1/2N}\right)$ such that $\widehat\nbf\sum_{i=1}^\infty i^{N-1}\chi(iq)\to 0$ and $h^{-d(1-\gamma)}\sum_{i=q}^\infty i^{N-1}(\chi(i))^{1-\gamma}\to 0$.
\end{description}

\begin{theo}\label{normquantile1}
Assume that (A1)-(A8), (C1) and (C2) hold. If there exists $c\geq 0$ such that $\widehat\nbf h^{d+4}\to c$, then
\begin{eqnarray*}
\sqrt{\widehat\nbf h^d}\left(\mu_{p,\nbf}(x)-\mu_p(x)\right)\stackrel{\mathcal{L}}{\to} \mathcal{N}\left(c\frac{B(x,\mu_p(x))}{f(\mu_p(x)|x)},\frac{\sigma^2(x,\mu_p(x))}{(f(\mu_p(x)|x))^2}\right),
\end{eqnarray*}
where
\begin{eqnarray}
\label{biais}
B(x,y) &= & \frac{1}{2}\Bigg\{\sum_{i,j=1}^d\left[\frac{\partial^2F(y|x)}{\partial x_i\partial x_j}+\frac{2}{g(x)}\frac{\partial g(x)}{\partial x_i}\frac{\partial F(y|x)}{\partial x_j}\right]\int_{\mathbb{R}^d}\|s\|^2K(s)ds \nonumber\\ \mbox{} & & +\frac{\partial^2 F(y|x)}{\partial y^2}\int_{\mathbb{R}}t^2w(t)dt\Bigg\}.\\
\sigma^2(x,y) & = & \frac{F(y|x)\left[1-F(y|x)\right]}{g(x)}\int_{\mathbb{R}^d}K^2(z)dz.
\label{variance}
\end{eqnarray}
\end{theo}
\subsection{Prediction}\label{sec_pred}
Let $(\xi_{\mathbf{i}},\,\mathbf{i}\in \Ibf)$ be a $\mathbb{R}-$valued strictly stationary random spatial process, assumed to be bounded, observable over a region
$\Ibf\subset\mathbb{N}^{N}$ and observed over a subset $\Obf$ of $\Ibf$, $\mathbf{n}=(n_{1},...,n_{N})\in\mathbb{N}^{N}$. The aim of this section is to predict
$\xi_{\mathbf{i}_{0}}$, at a given fixed point $\mathbf{i}_{0}$ not in $\Obf\subset\mathbb{N}^{N}$.
In practice (e.g. for simplicity), we expect that $\xi_{\mathbf{i}_{0}}$ depends only on the values of the process on a bounded neighborhood
$\mathcal{V}_{\mathbf{i}_{0}}\subset\Obf$. In other words, we expect that the process $(\xi_{\mathbf{i}})$ satisfies a Markov property, see for example Biau and Cadre
\cite{biau}, Dabo-Niang and Yao \cite{dabo}. Moreover, we assume that $\mathcal{V}_{\mathbf{i}_{0}}=\mathcal{V}+\mathbf{i}_{0}$, where
$\mathcal{V}$ is a fixed bounded set of sites that does not contain $0$. It is well known that the best predictor of $\xi_{\mathbf{i}_{0}}$ given the data in
$\mathcal{V}_{\mathbf{i}_{0}}$ in the sense of \emph{mean-square error} is
\[ E(\xi_{\mathbf{i}_{0}}|\xi_{\mathbf{i}},\mathbf{i}\in\mathcal{V}_{\mathbf{i}_{0}}).
\]

\noindent Let $\mathcal{V}_{\mathbf{i}}=\mathcal{V}+\mathbf{i}=\{\mathbf{u}+\mathbf{i},\,\mathbf{u}\in\mathcal{V}\}$ for each $\mathbf{i}\in\mathbb{N}^{N}$, and $d$ be the
cardinal of $\mathcal{V}$ ($d$ is also the cardinal of each $\mathcal{V}_{\mathbf{i}}$). To define a predictor of $\xi_{\mathbf{i}_{0}}$, let us consider the
$\mathbb{R}^d$-valued random variables $\tilde{\xi}_{\mathbf{i}}=\{\xi_{\mathbf{u}},\; {\mathbf{u}}\in \mathcal{V}_{\mathbf{i}}\subset\Obf\}$. The notation of the previous sections are used by setting $X_{\mathbf{i}}=\tilde{\xi}_{\mathbf{i}},\;Y_{\mathbf{i}}=\xi_{\mathbf{i}},\,\mathbf{i}\in\mathbb{N}^{N}$.\\ As a
predictor of $\xi_{\mathbf{i_0}}$, we take the conditional quantile estimate $\widehat{\xi}_{\mathbf{i_0}}=\mu_{p,\nbf}(\tilde{\xi}_{\mathbf{i_0}})$ of order $p$,
particularly the conditional median $p=0.5$.
%The quantile estimate $\mu_{p,\nbf}(\tilde{\xi}_{\mathbf{i_0}})$ is defined as the root of the equation
% \begin{eqnarray*}\label{double_estimator}
%F_{\nbf}(y|\tilde{\xi}_{\mathbf{i_0}})=p.
%\end{eqnarray*}
%\begin{eqnarray*}
%F_{\nbf}(y|\tilde{\xi}_{\mathbf{i_0}})  &= & \frac{\psi_{\nbf}(\tilde{\xi}_{\mathbf{i_0}},y)}{g_\nbf(\tilde{\xi}_{\mathbf{i_0}})}\1_{\left\{g_\nbf(\tilde{\xi}_{\mathbf{i_0}})\neq 0\right\}},\\
%\psi_{\nbf}(\tilde{\xi}_{\mathbf{i_0}},y) & = & \frac{1}{\widehat\nbf
%h^{d+1}}\sum_{\ibf\in\Ibf}K\left(\frac{\tilde{\xi}_{\mathbf{i_0}}-\tilde{\xi}_{\mathbf{i}}}{h}\right)\int_{-\infty}^yw\left(\frac{y-\tilde{\xi}_{\mathbf{i}}}{h}\right)dz.
%\end{eqnarray*}
%\begin{eqnarray*}
%g_\nbf(\tilde{\xi}_{\mathbf{i_0}}) & = & \frac{1}{\widehat\nbf h^d}\sum_{\ibf\in\Ibf}K\left(\frac{\tilde{\xi}_{\mathbf{i_0}}-\tilde{\xi}_{\mathbf{i}}}{h}\right),
%\end{eqnarray*}
We deduce from the previous consistency results, the following corollary that gives the convergence of the predictor
$\widehat{\xi}_{\mathbf{i_0}}$.
\begin{cor}
\label{prediction-markov}
\begin{description}
\item i) Under the conditions of Corollary \ref{convl1}, we have
\begin{eqnarray*}
\Ebf\left[\{\mu_{p,\nbf}(\tilde{\xi}_{\mathbf{i_0}})-\mu_{p}(\tilde{\xi}_{\mathbf{i_0}})\}\1_{\{\tilde{\xi}_{\mathbf{i_0}}\in S\}}\right]\stackrel{a.s.}{\to}
0.
\end{eqnarray*}
\item ii)  Under the conditions of Theorem \ref{normquantile1}, and if $\widehat{n}h^{d+4} \to 0$, then
\begin{eqnarray*}
\sqrt{\widehat\nbf h^d}\left(\mu_{p,\nbf}(\tilde{\xi}_{\mathbf{i_0}})-\mu_{p}(\tilde{\xi}_{\mathbf{i_0}})\right)\stackrel{\mathcal{L}}{\to} \mathcal{N}\left(%c\frac{B(\tilde{\xi}_{\mathbf{i_0}},\mu_{p}(\tilde{\xi}_{\mathbf{i_0}}))}{f(\mu_{p}(\tilde{\xi}_{\mathbf{i_0}})|\tilde{\xi}_{\mathbf{i_0}})},
0,\frac{\widehat{\sigma}^2(\tilde{\xi}_{\mathbf{i_0}},\mu_{p}(\tilde{\xi}_{\mathbf{i_0}}))}{(f(\mu_{p}(\tilde{\xi}_{\mathbf{i_0}})|\tilde{\xi}_{\mathbf{i_0}}))^2}\right),
\end{eqnarray*}
\end{description}
\end{cor}
where
\begin{eqnarray*}
\widehat{\sigma}^2(x,y)=\frac{F_{\nbf}(y|x)\left[1-F_{\nbf}(y|x)\right]}{g_{\nbf}(x)}\int_{\mathbb{R}^d}K^2(z)dz.
\end{eqnarray*}

These consistency results permit to have an approximation of an $1-\alpha$ confidence interval of $\xi_{\mathbf{i_0}}$ given by
$\widehat{I}_\alpha=[a_{-}(\tilde{\xi}_{\mathbf{i_0}}), a_{+}(\tilde{\xi}_{\mathbf{i_0}})]$, where
\begin{equation}
a_{\pm}(\tilde{\xi}_{\mathbf{i_0}})=\mu_{p,\nbf}(\tilde{\xi}_{\mathbf{i_0}})\pm
Q_{1-\frac{\alpha}{2}}\frac{\sigma(\tilde{\xi}_{\mathbf{i_0}},\mu_{p,\nbf}(\tilde{\xi}_{\mathbf{i_0}}))}{\sqrt{\widehat\nbf
h^d}f_\nbf(\mu_{p,\nbf}(\tilde{\xi}_{\mathbf{i_0}})|\tilde{\xi}_{\mathbf{i_0}})}, \label{intcon}
\end{equation}
where $Q_{\zeta}$ denotes the $\zeta$-quantile of the  standard normal distribution, and the unknown parameters (of the asymptotic variance in Corollary
\ref{prediction-markov}) are replaced by
kernel estimates.\\
Note also that the quantiles of order $p_1$ and $p_2$ ($p_1<p_2$) can be used to construct a predictive interval that consists of the $(p_2-p_1)100\%$ confidence
interval with bounds $\mu_{p_1,\nbf}(\tilde{\xi}_{\mathbf{i_0}})$ and $\mu_{p_2,\nbf}(\tilde{\xi}_{\mathbf{i_0}})$.

\section{A simulation study}
\begin{figure}[H]

 \centering
 \begin{tabular}{cc}
   &\includegraphics[scale=0.31]{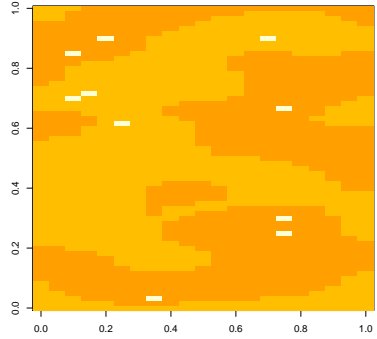}
   \end{tabular}
      \caption{The random field $(\xi_{\mathbf{i}})$ with non observed values $\xi_{\mathbf{i}_{1}},\ldots,\xi_{\mathbf{i}_{10}}$ in the white rectangular cases.}
       \label{champ1}
 \end{figure}

In this section, we study the performance of the conditional quantile predictor introduced in the previous section  towards some simulations. Let us denoted by $GRF(m,\,\sigma^{2},s)$ a Gaussian random field with mean
$m$ and covariance function defined by
$$\vartheta(h)=\sigma^{2}\exp\left\{-\left(\frac{\|h\|}{s}\right)^{2}\right\},h\in\mathbb{R}^{2}.$$
Set
\begin{eqnarray}
\label{rectangle}
\Ibf&=&\{\ibf=(i,j) \in (\mathbb{N^*})^2, \ \  1\leq i\leq 61, \ \ 1\leq j\leq 61 \}\\
 \mathcal{O}_{\mathbf{n}}&=&\{\ibf=(i,j) \in (\mathbb{N^*})^2, \ \  1\leq i\leq 21, \ \ 1\leq j\leq 21\} \cup \{\ibf=(22,j), \  1\leq j\leq 15\}.
 \label{observe}
 \end{eqnarray}
We consider a random field $(\xi_{\mathbf{i}})_{\ibf\in \Ibf}$ from the following model
 \begin{eqnarray}
\label{model} \xi_{\mathbf{i}}=U_{\mathbf{i}}*\left(\sin(2X_{\mathbf{i}})+2\exp\{-(16X_{\mathbf{i}})^2\}\right)+Z_{\mathbf{i}},\; \mathbf{i}\in \mathbb{N}^2
\end{eqnarray}
where  $X=(X_{\mathbf{i}})_{\ibf\in \Ibf}$ is a $GRF(0,\,5,\ 3)$, $Z=(Z_{\mathbf{i}})_{\ibf\in \Ibf}$ is a $GRF(0,\,0.1,\ 5)$ independent of $X$ and
$U_{\mathbf{i}}=\frac{1}{\hat{\mathbf{n}}}\sum_{\mathbf{j}\in \Ibf}\exp\left(-\frac{\left\Vert \mathbf{i}-\mathbf{j}\right\Vert }{2}\right)$. The choice of $U_\ibf$
in the model (\ref{model}) is motivated by a reinforcement of the spatial local dependency. The field $(\xi_{\ibf},\ibf \in
\Ibf)$ is observable over the rectangular region $\Ibf$ and observed over the subset $ \mathcal{O}_{\mathbf{n}}$ defined in \eqref{rectangle} and \eqref{observe}.\\
We want to predict the values $\xi_{\mathbf{i}_{1}},\ldots,\xi_{\mathbf{i}_{m}}$ at given fixed sites $\mathbf{i}_{1}, \ldots, \mathbf{i}_{m}$ not in
$\mathcal{O}_{\mathbf{n}}$, with $m=10$. The sample obtained from model \eqref{model}, observed in $\mathcal{O}_{\mathbf{n}}$ is plotted in Figure \ref{champ1} below with
the $10$ non observable values of the field at $\mathbf{i}_{1}, \ldots, \mathbf{i}_{m}$.
As explained in Section \ref{sec_pred}, for any $k\in\{1,\ldots,m\}$, we
take the conditional quantile estimate $\widehat{\xi}_{\mathbf{i}_k}=\mu_{p,\nbf}(\tilde\xi_{\mathbf{i}_k})$ as a predictor of ${\xi}_{\mathbf{i}_k}$, where
$\tilde\xi_{\mathbf{i}_k}$ are observed on $\Obf$ and the vicinity $\mathcal{V}=\{-1,1\}\times \{-1,1\}$ or $\{-2,-1,1,2\}\times\{-2,-1,1,2\}$.\\ To compute
$\mu_{p,\nbf}$, we select the standard normal density as kernel $K$ and the Epanechnikov kernel as $w$. For the bandwidth selection, we use the rule developed in Yu
and Jones \cite{yu},
\begin{eqnarray*}
h_\nbf=h_{\mbox{mean}}\left(\frac{p(1-p)}{\phi(\Phi^{-1}(p))^2}\right)^{1/5},
\end{eqnarray*}
where $h_{\mbox{mean}}$ is the bandwidth for kernel smoothing estimation of the regression mean,
$\phi$ and $\Phi$, are respectively, the standard normal density and distribution function.\\
 To evaluate the
performance of the predictor $\widehat{\xi}_{\mathbf{i}_k}$, we compute the mean absolute error (MAE):
\begin{eqnarray*}
MAE = \frac{1}{m}\sum_{k=1}^m|\widehat{\xi}_{\mathbf{i}_k}-\xi_{\mathbf{i}_{k}}|.
\end{eqnarray*}
 The following Table gives the predictors of $\xi_{\mathbf{i}_{k}},\; k=1,...,m$ for $p\in \{0.05,\ 0.5,\ 0.95\}$,
 $\mathcal{V}=\{-1,1\}\times \{-1,1\}$ on the left, $\mathcal{V}=\{-2,-1,1,2\}\times\{-2,-1,1,2\}$ on the right and the prediction error.
 \begin{table}[H]
 \caption{Predictive data for $\mathcal{V}=\{-1,1\}\times\{-1,1\}$ on the left and $\mathcal{V}=\{-2,-1,1,2\}\times\{-2,-1,1,2\}$ on the right.}
 \centering
 \begin{tabular}{c|cccc||cccc}
\hline
&$p=0.05$& True data &$p=0.5$& $p=0.95$&$p=0.05$& True data &$p=0.5$& $p=0.95$\\\hline
&0.1653 &$\mathbf{0.2009}$   &0.1930  &  0.2192  &0.1835 &$\mathbf{0.2009}$   &0.2129 &  0.2362   \\
&-0.2553&$\mathbf{ -0.2315}$  & -0.2289 &  -0.1862&-0.2195 &$\mathbf{-0.2315}$&-0.1984 & -0.1766  \\
& 0.1516&$\mathbf{0.1966}$   &0.1990& 0.2362  &0.1912&$\mathbf{0.1966}$ &0.2129&  0.2362   \\
& -0.5313&$\mathbf{ -0.4906}$   & -0.5062 & -0.4782 &-0.5472& $\mathbf{-0.4906}$  &-0.5313  &-0.5033   \\
& 0.2693&$\mathbf{  0.2901}$    & 0.2929 &  0.3168 & 0.2237&$\mathbf{0.2901}$    & 0.2465& 0.2676   \\
&-0.2748 & $\mathbf{ -0.2535}$ &-0.2527  & -0.2289 &-0.2838&$\mathbf{-0.2535}$  &  -0.2606 & -0.2401   \\
&0.3696&$\mathbf{0.3941}$  &0.4007 & 0.4269  &0.3805&$\mathbf{0.3941}$  &0.3834& 0.4269   \\
&-0.5539&$\mathbf{-0.5177}$   &-0.5295   & -0.5062 &-0.5472&$\mathbf{-0.5177}$ & -0.5313&-0.5033   \\
 &-0.3678&$\mathbf{-0.3217}$  & -0.3463  & -0.3193 &-0.3637&$\mathbf{-0.3217}$  &-0.3487  & -0.3231   \\
& -0.2983&$\mathbf{-0.2843}$ &-0.2702    &  -0.2455&-0.3096& $\mathbf{-0.2843}$ &-0.2863 &-0.2671\\\hline\hline
 MAE & 0.0308 & &0.0089 & 0.0252&   0.0298  & &0.0206   & 0.0244
 \end{tabular}
 \end{table}
% \begin{table}[H]
 %  \footnotesize{
 % \caption{Predictive data for $\widehat\nbf=36\times36$ and $\mathcal{V}=\{-2,-1,1,2\}\times\{-2,-1,1,2\}$.}
 %\centering
 %\begin{tabular}{c|cccc}\hline
%&p=0.05& True data & p=0.5& p=0.95\\\hline
%&0.1835 &$\mathbf{0.2009}$   &0.2129 &  0.2362   \\
%&-0.2195 &$\mathbf{-0.2315}$&-0.1984 & -0.1766  \\
%&0.1912&$\mathbf{0.1966}$ &0.2129&  0.2362   \\
%&-0.5472& $\mathbf{-0.4906}$  &-0.5313  &-0.5033   \\
%& 0.2237&$\mathbf{0.2901}$    & 0.2465& 0.2676   \\
%&-0.2838&$\mathbf{-0.2535}$  &  -0.2606 & -0.2401   \\
%&0.3805&$\mathbf{0.3941}$  &0.3834& 0.4269   \\
%&-0.5472&$\mathbf{-0.5177}$ & -0.5313&-0.5033   \\
%&-0.3637&$\mathbf{-0.3217}$  &-0.3487  & -0.3231   \\
%&-0.3096& $\mathbf{-0.2843}$ &-0.2863 &-0.2671\\\hline\hline
%MAE &   0.0298  & &0.0206   & 0.0244
% \end{tabular}}
 %\end{table}
We derive from the results of Tables 1 a $90\%$ predictive interval where the extremities are the $5\%$ and $95\%$ quantiles estimates, for each of the $10$
prediction sites (see Section \ref{sec_pred} for more details). Note that these $90\%$ predictive intervals contain the true values. The average length for the $10$ intervals is
$ 0.05604$ for  $\mathcal{V}=\{-1,1\}\times \{-1,1\}$ and $0.04455$ for $\mathcal{V}=\{-2,-1,1,2\}\times\{-2,-1,1,2\}$.\\
The numerical results show that our proposed predictor gives good results on the above simulated field. A next step would be to apply the predictor to a spatial real data and deserves futur investigations.

\section{Appendix}
The letter $C$ will be used to denote constants whose values are
unimportant. Before proving the main results, let us give the following notations:
\begin{eqnarray*}
\psi(x,y)=\int_{-\infty}^yf(x,z)dz, \ \ K_1(y)=\int_{-\infty}^yw(z)dz.
\end{eqnarray*}
Next, we need the following results.
%\begin{lem}\label{convmargas} (Carbon {\it et. al} (1997))
%Assume that Assumptions (A1), (A2), (A3)i), (A5) and (A7) hold.
%\begin{description}
%\item (i) If (A9) and (\ref{mixing1}) are satisfied,
%\item (ii) or if (A10) and (\ref{mixing2}) are satisfied
%\end{description}
%then
%\begin{eqnarray*}
%\sup_{x\in S}|g_\nbf(x)-g(x)|=O\left(\Omega_\nbf\right)\ \ \mbox{a.s.}
%\end{eqnarray*}
%\end{lem}

\begin{lem}\label{psibias}
Assume that (A3)-(A7) hold, then,
\begin{eqnarray*}
\sup_{y \in \mathcal{V}}\sup_{x\in
S}\left|\Ebf\psi_{\nbf}(x,y)-\psi(x,y)\right| &= &
O\left(\Omega_\nbf\right).
\end{eqnarray*}
\end{lem}
The proof of Lemma \ref{psibias} is classical (see for example Matzner and L\o ber \cite{matzner}).
%A simple computation shows that Assumption (A7) implies that $\widehat\nbf h^{d+4}(\log\nbf)^{-1}=o(1)$ so that $h^2=o\left(\Omega_\nbf\right)$.
%Note that,
%\begin{eqnarray*}
%\left|\Ebf\psi_\nbf(x,y)-\psi(x,y)\right| & = & \left|\frac{1}{\widehat\nbf h^d}\sum_{\ibf \in \Ibf}\int_{\mathbb{R}^{d+1}}K\left(\frac{x-z}{h}\right)K_1\left(\frac{y-v}{h}\right)f(z,v)dzdv-\int_{-\infty}^yf(x,u)du\right|\\&= &
% \left|\int_{-\infty}^y\int_{\mathbb{R}^{d+1}}K(s)w(s')\left[f(x-sh,u-s'h)-f(x,u)\right]dudsds' \right|\\
% & \leq &\int_{]-\infty,y]\times\mathbb{R}^{d+1}} \left|K(s)w(s')\right|\left|f(x-sh,u-s'h)-f(x,u)\right|dudsds'\\ & \leq & Ch^2= O\left(\Omega_\nbf\right),
%\end{eqnarray*}
%by two order Taylor's expansion.

\begin{lem}\label{psivar2}
Let Assumptions (A1)-(A6) hold. If (A9) and (\ref{mixing1}) or (A10) and (\ref{mixing2}) are satisfied,
then,
\begin{eqnarray*}
\sup_{y \in \mathcal{V}}\sup_{x \in
S}\left|\psi_{\nbf}(x,y)-\Ebf\psi_{\nbf}(x,y)\right| &= & O\left(\Omega_\nbf\right)\ \ \mbox{a.s.}
\end{eqnarray*}

\end{lem}
\noindent\textbf{Proof of Lemma \ref{psivar2}.}
Define
\begin{eqnarray*}
Z_\ibf(x,y) & = & \frac{1}{\widehat \nbf h^{d+1}}K\left(\frac{x-X_\ibf}{h}\right)\int_{-\infty}^yw\left(\frac{t-Y_\ibf}{h}\right)dt\; =
 \frac{1}{\widehat \nbf h^d}K\left(\frac{x-X_\ibf}{h}\right)K_1\left(\frac{y-Y_\ibf}{h}\right),\\
 \Delta_\ibf(x,y) & = & Z_\ibf(x,y)-\Ebf Z_\ibf(x,y),
 \end{eqnarray*}
 then
$\displaystyle
S_\nbf(x,y)= \psi_{\nbf}(x,y)-\Ebf\psi_{\nbf}(x,y)=\sum_{\ibf\in \Ibf} \Delta_\ibf(x,y)$.
%Let
%$\displaystyle \tilde I_\nbf(x,y)  =  \sum_{\ibf\in \Ibf}\Ebf(\Delta_\ibf(x,y))^2  \ \ \mbox{and} \ \  \tilde R_\nbf(x,y)=\sum_{\underset{\ibf,\jbf\in\Ibf}{\ibf \neq
%\jbf}}\Ebf\Delta_\ibf(x,y) \Delta_\jbf(x,y). $

\begin{lem}\label{lem1}
Under Assumptions (A1), (A2), (A5) and (A6), we have
\begin{eqnarray*}
 \sum_{\ibf\in \Ibf}\Ebf(\Delta_\ibf(x,y))^2 +\sum_{\underset{\ibf,\jbf\in\Ibf}{\ibf \neq
\jbf}}\Ebf\Delta_\ibf(x,y) \Delta_\jbf(x,y)=O\left(\frac{1}{\widehat\nbf
h^d}\right).
\end{eqnarray*}
\end{lem}
The proof of Lemma \ref{lem1} is similar to the proof of Lemma 2.2 in Tran \cite{tran}.\\
%Now, define $\displaystyle S_\nbf(x,y)=\psi_{\nbf}(x,y)-\Ebf\psi_{\nbf}(x,y)$.\\

Let us introduce a spatial block decomposition that has been used by Tran \cite{tran} and Carbon et al. \cite{carbon}. Without loss of generality, assume that $n_i=2pq_i$ for $1\leq i\leq N$. The random variables $\Delta_\ibf(x,y)$ can be grouped into $2^Nq_1\times q_2\times\cdots\times q_N$ cubic blocks of side $p$. Denote
\begin{eqnarray*}
U^x(1,\nbf,\jbf,y) & = & \sum_{\underset{k=1,\ldots,N}{i_k=2j_kp+1}}^{(2j_k+1)p}\Delta_\ibf(x,y),\;
U^x(2,\nbf,\jbf,y)  =  \sum_{\underset{k=1,\ldots,N-1}{i_k=2j_kp+1}}^{(2j_k+1)p} \sum_{i_N=(2j_N+1)p+1}^{2(j_N+1)p}\Delta_\ibf(x,y),\\
U^x(3,\nbf,\jbf,y) & = & \sum_{\underset{k=1,\ldots,N-2}{i_k=2j_kp+1}}^{(2j_k+1)p}\sum_{i_{N-1}=(2j_{N-1}+1)p+1}^{2(j_{N-1}+1)p} \sum_{i_N=2j_Np+1}^{(2j_N+1)p}\Delta_\ibf(x,y),\\
U^x(4,\nbf,\jbf,y) & = & \sum_{\underset{k=1,\ldots,N-2}{i_k=2j_kp+1}}^{(2j_k+1)p} \sum_{i_{N-1}=(2j_{N-1}+1)p+1}^{2(j_{N-1}+1)p}\sum_{i_N=(2j_N+1)p+1}^{2(j_N+1)p}\Delta_\ibf(x,y),
\end{eqnarray*}
and so on. Finally, note that
\begin{eqnarray*}
U^x(2^{N-1},\nbf,\jbf,y) & = & \sum_{i_k=2j_kp+1}^{(2j_k+1)p} \sum_{i_N=2j_Np+1}^{(2j_N+1)p}\Delta_\ibf(x,y),\;\;U^x(2^N,\nbf,\jbf,y) =
\sum_{\underset{k=1,\ldots,N}{i_k=2j_kp+1}}^{(2j_k+1)p}\Delta_\ibf(x,y).
\end{eqnarray*}

For each integer $1\leq i\leq 2^N$, let
\begin{eqnarray*}
T^x(\nbf,i,y)=\sum_{\underset{k=1,\ldots,N}{j_k=0}}^{q_k-1}U^x(i,\nbf,\jbf,y),\; S_\nbf(x,y)=\sum_{i=1}^{2^N}T^x(\nbf,i,y).
\end{eqnarray*}
Observe that, for any $\varepsilon>0$
\begin{eqnarray}
\label{snd} \Pbf\left(|S_\nbf(x,y)|>\varepsilon\right) &= & \Pbf\left(\left|\sum_{i=1}^{2^N}T^x(\nbf,i,y)\right|>\varepsilon\right)\nonumber
  \leq
 2^N\Pbf\left(\left|T^x(\nbf,i,y)\right|>\varepsilon/2^N\right).
\end{eqnarray}
% A VOIR
Without loss of generality, we consider just the case where $i=1$ and we enumerate in an arbitrary way the $\widehat q=q_1\ldots q_N$ terms $U^x(1,\nbf,\jbf,y)$ of the sum $T^x(\nbf,1,y)$ that we call $W_1,\ldots,W_{\widehat q}$.  Note that $U^x(i,\nbf,\jbf,y)$ is measurable with respect to the $\sigma$-field generated by $Z_\ibf$, with $\ibf$ such that $2j_kp+1\leq i_k\leq (2j_k+1)p$, $k=1,\ldots, N$.\\
These sets of sites are separeted by a distance at least $p$ and since the $Z_\ibf$ are bounded, then we have for all $i=1,\ldots,\widehat q$, $\displaystyle
|W_i|\leq C(\widehat\nbf h^d)^{-1}p^N\|K\|_{\infty}.
$
Lemma 4.4 in Carbon et al. \cite{carbon} ensures that there exist independent random variables $W_1^*,\ldots,W_{\widehat q}^*$ such that for all $i=1,\ldots,\widehat q$,
\begin{eqnarray*}
\Ebf|W_i-W_i^*|\leq C(\widehat\nbf
h^d)^{-1}p^N\|K\|_{\infty}\phi(\widehat\nbf,p^N)\chi(p).
\end{eqnarray*}
Markov's inequality leads to
\begin{eqnarray}
\label{markov}
\Pbf\left(\sum_{i=1}^{\widehat
q}|W_i-W_i^*|>\varepsilon/2^{N+1}\right) \leq C2^{N+1}(\widehat\nbf
h^d)^{-1}p^N\widehat
q\|K\|_{\infty}\phi(\widehat\nbf,p^N)\varepsilon^{-1}\chi(p).
\end{eqnarray}
By Bernstein's inequality, we
have
\begin{eqnarray}
\label{bernstein}
\Pbf\left(|\sum_{i=1}^{\widehat q}W_i^*|>\varepsilon/2^{N+1}\right) \leq
2\exp\left\{\frac{-\varepsilon^2/(2^{N+1})^2}{4\sum_{i=1}^{\widehat q}\Ebf W_i^{*2}+2C(\widehat\nbf h^d)^{-1}p^N\|K\|_{\infty}\varepsilon/2^{N+1}}\right\}.
\end{eqnarray}
Combining \eqref{markov} and \eqref{bernstein}, we have
\begin{eqnarray*}
\Pbf\left(|S_\nbf(x,y)|>\varepsilon\right)
&\leq  &
2^N\Pbf\left(\sum_{i=1}^{\widehat q}|W_i-W_i^*|>\varepsilon/2^{N+1}\right)+2^N\Pbf\left(|\sum_{i=1}^{\widehat q}W_i^*|>\varepsilon/2^{N+1}\right)\\
& \leq &  2^{N+1}\exp\left\{\frac{-\varepsilon^2/(2^{N+1})^2}{4\sum_{i=1}^{\widehat q}\Ebf W_i^{*2}+2C(\widehat \nbf h^d)^{-1}p^N\|K\|_{\infty}\varepsilon/2^{N+1}}\right\} \\
\mbox{} && + C2^{2N+1}\phi(\widehat\nbf,p^N)(\widehat\nbf
h^d)^{-1}p^N\widehat q\|K\|_{\infty}\varepsilon^{-1}\chi(p).
\end{eqnarray*}
By Lemma \ref{lem1}, one has $\sum_{i=1}^{\widehat q}\Ebf W_i^{*2}=O(\frac{1}{\widehat\nbf h^d})$ and since $\widehat \nbf=2^Np^N\widehat q$, we have
\begin{eqnarray*}
\Pbf\left(|S_\nbf(x,y)|>\varepsilon\right) & \leq &
2^{N+1}\exp\left\{\frac{-\varepsilon^2\widehat \nbf
h^d}{2^{2N+4}C+2^{N+2}Cp^N\varepsilon}\right\} +
C2^{N+1}\phi(\widehat\nbf,p^N)h^{-d}\varepsilon^{-1}\chi(p).
\end{eqnarray*}
Let $\lambda >0$ and set $ \varepsilon =  \lambda\Omega_\nbf,\, p =  \Omega_\nbf^{-1/N}.$  For the first part of Lemma \ref{psivar2}, a  simple computation
shows that for sufficiently large $\widehat\nbf$,
\begin{eqnarray}
\label{vac1} \Pbf\left(|S_\nbf(x,y)|>\lambda\Omega_\nbf\right)& \leq &
2^{N+1}\exp\left\{\frac{-\lambda^2\log\widehat\nbf}{2^{2N+4}C+2^{N+2}C\lambda }\right\} + C2^{N+1}p^N h^{-d}\lambda^{-1}\Omega_\nbf^{-1}\chi(p) \nonumber \\ & \leq &
 %C\widehat\nbf^{-b}+C2^{N+1}\Omega_\nbf^{-1} h^{-d}\lambda^{-1}\Omega_\nbf^{-1}\Omega_\nbf^{\beta/N}\nonumber\\
 C\widehat\nbf^{-b}+C\lambda^{-1}h^{-d}\Omega_\nbf^{\frac{\beta-2N}{N}},
\end{eqnarray}
with $b>0$.
Analogously, for the second part, as above we have
\begin{eqnarray}
\label{vac2} \Pbf\left(|S_\nbf(x,y)|>\lambda\Omega_\nbf\right)& \leq &
2^{N+1}\exp\left\{\frac{-\lambda^2\log
\widehat\nbf}{2^{2N+4}C+2^{N+2}C\lambda }\right\} +C2^{N+1}\widehat\nbf^{\kappa}h^{-d}\lambda^{-1}\Omega_\nbf^{-1}\chi(p)\nonumber\\
& \leq & %C\widehat\nbf^{-b}+C2^{N+1}\widehat\nbf^{\kappa}h^{-d}\lambda^{-1}\Omega_\nbf^{-1}\Omega_\nbf^{\beta/N}\nonumber \\ & \leq &
C\widehat\nbf^{-b}+C\lambda^{-1}\widehat\nbf^{\kappa}h^{-d}\Omega_\nbf^{\frac{\beta-N}{N}}.
\end{eqnarray}
Now, set $R_\nbf=h^{d+1}\Omega_\nbf$ and $r_\nbf=\left(\frac{h^{d+2}}{\nbf}\right)^{1/2}$. Since $S$ is a compact, it can be covered with $d_\nbf$
cubes $B_k$ having sides of length $R_\nbf$ and center at $x_k$ with $d_\nbf\leq CR_\nbf^{-d}$. The compact set $\mathcal{V}$ can be covered with $l_\nbf$ intervalls $I_l$ having length $r_\nbf$ and center at $y_l$, with $l_\nbf \leq Cr_\nbf^{-1}$.
We have
\begin{eqnarray}
\label{dn}
d_\nbf l_\nbf\leq C\nbf^{1/2}h^{-\frac{d+2}{2}}\left(h^{d+1}\Omega_\nbf\right)^{-d-1}.
\end{eqnarray}
Define $S_{1\nbf} = \sup_{y\in \mathcal{V}} \max_k\sup_{x\in B_k}\left|\psi_{\nbf}(x,y)-\psi_{\nbf}(x_k,y)\right|$
\begin{eqnarray*}
S_{2\nbf} &=&  \max_k\max_l\sup_{y\in I_l}\left|\psi_{\nbf}(x_k,y)-\psi_{\nbf}(x_k,y_l)\right|,\; S_{3\nbf} =\max_k\max_l\sup_{y\in
I_l}\left|\Ebf\psi_{\nbf}(x_k,y_l)-\Ebf\psi_{\nbf}(x_k,y)\right|\\
S_{4\nbf} &=&  \sup_{y \in \mathcal{V}}\max_k\sup_{x\in
B_k}\left|\Ebf\psi_{\nbf}(x_k,y)-\Ebf\psi_{\nbf}(x,y)\right|,\; S_{5\nbf} =\max_k\max_l\left|\psi_{\nbf}(x_k,y_l)-\Ebf\psi_{\nbf}(x_k,y_l)\right|.
\end{eqnarray*}
Then, we can write $\displaystyle\sup_{y \in \mathcal{V}}\sup_{x\in S}\left|\psi_{\nbf}(x,y)-\Ebf\psi_{\nbf}(x,y)\right| \leq S_{1\nbf}+S_{2\nbf}+S_{3\nbf}+S_{4\nbf}+S_{5\nbf}$.
The proof of Lemma \ref{psivar2} follows easily from the combination of the two following lemmas.
\begin{lem}\label{sn1sn2}
Under Assumptions (A5) and (A6),
\begin{eqnarray*}
S_{i\nbf} & = & O\left(\Omega_\nbf\right) \ \  \mbox{a.s.},\; i=1,2,3,4.\\
\end{eqnarray*}
\end{lem}
\begin{lem}\label{sn3}
Assume that (A1), (A2), (A5) and (A6) hold. If (A9) and (\ref{mixing1}) or A(10) and (\ref{mixing2}) are satisfied,
then
\begin{eqnarray*}
 S_{5\nbf} & = & O\left(\Omega_\nbf\right) \ \  \mbox{a.s.}
\end{eqnarray*}
\end{lem}
\noindent\textbf{Proof of Lemma \ref{sn1sn2}.}
On one hand, since the kernel $K$ satisfies
the Lipschitz condition, we have clearly
\begin{eqnarray*}
\left|\psi_{\nbf}(x,y)-\psi_{\nbf}(x_k,y)\right|& \leq & C\widehat\nbf^{-1} h^{-d-1}\sum_{\ibf\in\Ibf}\|x-x_k\| \leq  Ch^{-d-1}R_\nbf=O\left(\Omega_\nbf\right).
\end{eqnarray*}
On the other hand, observe that $h^{-d-1}r_\nbf=\left(\widehat\nbf h^d\right)^{-1/2}=O\left(\Omega_\nbf\right)$. Since $w$ satisfies the Lipschitz condition,
\begin{eqnarray*}
\left|\psi_{\nbf}(x_k,y)-\psi_{\nbf}(x_k,y_l)\right|& \leq & C\widehat\nbf^{-1} h^{-d-1}\sum_{\ibf\in\Ibf}\left|y-y_l\right|
  \leq  Ch^{-d-1}r_\nbf=O\left(\Omega_\nbf\right),
\end{eqnarray*}
which gives the proof of Lemma \ref{sn1sn2}.\\

\noindent\textbf{Proof of Lemma \ref{sn3}.}
For $\varepsilon >0$, we have
\begin{eqnarray*}
\Pbf\left(\max_l\max_k\left|\psi_{\nbf}(x_k,y_l)-\Ebf\psi_{\nbf}(x_k,y_l)\right|>\varepsilon\right) &\leq &l_\nbf d_\nbf\Pbf\left(\left|\psi_{\nbf}(x_k,y_l)-\Ebf\psi_{\nbf}(x_k,y_l)\right|>\varepsilon\right).
\end{eqnarray*}
Setting $\varepsilon =\lambda\Omega_\nbf$ with $\lambda>0$ and taking into account \eqref{vac1} and \eqref{vac2}, it suffices to show that
$l_\nbf d_\nbf\widehat\nbf^{-b}\widehat\nbf u(\nbf)\to 0$ and $l_\nbf d_\nbf
h^{-d}\Omega_\nbf^{\frac{\beta-2N}{N}}\widehat\nbf u(\nbf)\to 0$, or $l_\nbf d_\nbf \widehat\nbf^{\kappa}h^{-d}\Omega_\nbf^{\frac{\beta-N}{N}}\widehat\nbf u(\nbf)\to 0$. \\
First, observe that condition $\widehat\nbf h^d\to \infty$ implies that
$\widehat\nbf>Ch^{-d}$, so that
\begin{eqnarray}
\label{eq}
\widehat\nbf^{(d+1)(d+2)/2d}>Ch^{-(d+1)(d+2)/2}.
\end{eqnarray}
 Using \eqref{dn} and \eqref{eq}, we have
\begin{eqnarray*}
l_\nbf d_\nbf \widehat\nbf^{-b}\widehat\nbf u(\nbf) &\leq& Ch^{-(d+1)\left(d+2\right)/2}(\log\widehat\nbf)^{-d/2}\widehat\nbf^{1/2}\widehat\nbf^{-b}\widehat\nbf^{d/2}\widehat\nbf u(\nbf)\\& \leq &
C\widehat\nbf^{(d^2+3d+1)/d-b}(\log\widehat\nbf)^{-d/2} u(\nbf),
\end{eqnarray*}
which goes to $0$ if $b>(d^2+3d+1)/d$. \\
Next, again \eqref{dn} and a computation show that
\begin{eqnarray*}
l_\nbf d_\nbf
h^{-d}\Omega_\nbf^{\frac{\beta-2N}{N}}\widehat\nbf u(\nbf)
 &\leq & C\left[\widehat\nbf(\log\widehat\nbf)^{\frac{\beta-N(d+2)}{N(d+5)-\beta}}u(\nbf)^{\frac{2N}{N(d+5)-\beta}}h^{\frac{-N(d+1)(d+2)-d\beta)}{N(d+5)-\beta}}\right]^{\frac{N(d+5)-\beta}{2N}},
\end{eqnarray*}
which goes to $0$ by Assumption (A9) and $\beta>N(d+5)$.\\
Analogously, \eqref{dn} and a computation show that
\begin{eqnarray*}
l_\nbf d_\nbf \widehat\nbf^{\kappa}h^{-d}\Omega_\nbf^{\frac{\beta-N}{N}}\widehat\nbf u(\nbf)
 &\leq & C\left[\widehat\nbf(\log\widehat\nbf)^{\frac{\beta-N(d+1)}{N(d+4+2\kappa)-\beta}}u(\nbf)^{\frac{2N}{N(d+4+2\kappa)-\beta}}h^{\frac{-N(d^2+4d+2)-d\beta)}{N(d+4+2\kappa)-\beta}}\right]
 ^{\frac{N(d+4+2\kappa)-\beta}{2N}},
\end{eqnarray*}
which goes to $0$ by Assumption (A10) and $\beta>N(d+4+2\kappa)$.
The conclusion of Lemma \ref{sn3} follows from the Borel Cantelli's Lemma.\\

\noindent \textbf{Proof of Theorems \ref{convcondps} and \ref{convmu}.}
First, from Carbon et al. \cite{carbon}, we have
$
\sup_{x\in S}|g_\nbf(x)-g(x)|=O\left(\Omega_\nbf\right)\  \mbox{a.s.}
$
Now, a standard decomposition gives
\begin{eqnarray*}
\label{decomp}
\lefteqn{\sup_{y \in \mathcal{V}}\sup_{x \in S}|F_{\nbf}(y|x)-F(y|x)|}\\& \leq  &
\frac{1}{\inf_{x\in S}g_\nbf(x)}\left\{\sup_{y \in \mathcal{V}}\sup_{x \in S}|\psi_{\nbf}(x,y)-\psi(x,y)|+\sup_{y \in \mathcal{V}}\sup_{x \in S}\left(F(y|x)|g(x)-g_\nbf(x)|\right)\right\}\nonumber\\
& \leq &
\frac{1}{\inf_{x \in S}g_\nbf(x)}\Bigg\{\sup_{y \in \mathcal{V}}\sup_{x \in S}|\psi_{\nbf}(x,y)-\Ebf\psi_{\nbf}(x,y)|+\sup_{y \in \mathcal{V}}\sup_{x \in S}|\Ebf\psi_{\nbf}(x,y)-\psi(x,y)|
+\sup_{x \in S}|g(x)-g_\nbf(x)|\Bigg\}.
\end{eqnarray*}
Since by (A1), $g_\nbf(x)$ is bounded away from $0$, Theorem \ref{convcondps} follows from the preceding inequality, Lemmas \ref{psibias} and \ref{psivar2}.
%\subsection*{Proof of Theorem \ref{convmu} }
Next, from (A8), to prove Theorem \ref{convmu}, it suffices to show that $\sup_{x\in S}\left|F(\mu_{p,\nbf}(x)|x)-F(\mu_p(x)|x)\right|\to 0$.
We have
\begin{eqnarray*}
\left|F(\mu_{p,\nbf}(x)|x)-F(\mu_p(x)|x)\right| &\leq &\left|F(\mu_{p,\nbf}(x)|x)-F_{\nbf}(\mu_{p,\nbf}(x)|x)\right|+\left|F_{\nbf}(\mu_{p,\nbf}(x)|x)-F(\mu_p(x)|x)\right| \\
& \leq & \left|F(\mu_{p,\nbf}(x)|x)-F_{\nbf}(\mu_{p,\nbf}(x)|x)\right|  \leq  \sup_{y \in \mathcal{V}}\left|F_{\nbf}(y|x)-F(y|x)\right|.
\end{eqnarray*}
Thus,
$
\sup_{x\in S}\left|F(\mu_{p,\nbf}(x)|x)-F(\mu_p(x)|x)\right|
 \leq  \sup_{x\in S} \sup_{y \in \mathcal{V}}\left|F_{\nbf}(y|x)-F(y|x)\right|,
$
so that Theorem \ref{convmu}  follows from an application of Theorem \ref{convcondps}.\\
% in conjonction with the inequality
%\begin{eqnarray*}
%\sum_\nbf \Pbf\left\{\sup_{x \in S}|\mu_{p,\nbf}(x)-\mu_p(x)|\geq \varepsilon\right\} \leq \sum_\nbf\Pbf\left\{\sup_{x \in S}\sup_{y \in %\mathcal{V}}|F_{\nbf}(y|x)-F(y|x)|\geq \varepsilon\right\}.
%\end{eqnarray*}
%%%%%%%%%%%%%%%%%estimateur 2%%%%%%%%%%%%%%%%%%

\noindent\textbf{Proof of Corollary \ref{convl1}.}
First, by Lyapounov's inequality, we have
\begin{eqnarray*}
\Ebf\left[\mu_{p,\nbf}(X_\nbf)-\mu_p(X_\nbf)\right] & \leq & \Ebf\left[|\mu_{p,\nbf}(X_\nbf)-\mu_p(X_\nbf)|\right]   \leq
\left(\Ebf\left[(\mu_{p,\nbf}(X_\nbf)-\mu_p(X_\nbf))^2\right] \right)^{1/2},
\end{eqnarray*}
so, we can write
\begin{eqnarray*}
\Ebf\left[\left(\mu_{p,\nbf}(X_\nbf)-\mu_p(X_\nbf)\right)\1_{\{X_\nbf \in
S\}}\right] & \leq &
\left(\Ebf\left[(\mu_{p,\nbf}(X_\nbf)-\mu_p(X_\nbf))^2\1_{\{X_\nbf \in
S\}}\right] \right)^{1/2}.
\end{eqnarray*}
Then we also have
\begin{eqnarray*}
\Ebf\left[(\mu_{p,\nbf}(X_\nbf)-\mu_p(X_\nbf))^2\1_{\{X_\nbf \in
S\}}\right] & \leq & \Ebf\left[\sup_{x\in
S}\left(\mu_{p,\nbf}(x)-\mu_p(x)\right)^2\right].
\end{eqnarray*}
An integration by parts gives
\begin{eqnarray*}
 \Ebf\left[\sup_{x\in S}\left(\mu_{p,\nbf}(x)-\mu_p(x)\right)^2\right] & \leq & 2\int_{0}^{\infty} v\Pbf\left(\sup_{x\in S}\left|\mu_{p,\nbf}(x)-\mu_p(x)\right|>v\right)dv.
\end{eqnarray*}
Using Assumption (A8), we have that for $\nbf$ large enough,
\begin{eqnarray*}
\int_{0}^{\infty} v\Pbf\left(\sup_{x\in S}\left|\mu_{p,\nbf}(x)-\mu_p(x)\right|>v\right)dv &\leq& 2\int_{0}^{\infty} t\Pbf\left(\sup_{x\in S}\left|F(\mu_p(x)|x)-F_{\nbf}(\mu_p(x)|x)\right|>t\right)dt\\ & \leq &
2\int_{0}^{\infty} t\Pbf\left(\sup_{x\in S}\sup_{y \in \mathcal{V}}\left|F(y)|x)-F_{\nbf}(y|x)\right|>t\right)dt,%\\ & =& 2\Omega_{\nbf}^2\int_{0}^{\infty} \lambda\Pbf\left(\sup_{x\in S}\sup_{y \in \mathcal{V}}\left|F(y)|x)-F_{\nbf}(y|x)\right|>\lambda\Omega_{\nbf}\right)d\lambda,
\end{eqnarray*}
and Corollary \ref{convl1} follows from Theorem \ref{convcondps}.\\

\noindent\textbf{Proof of Theorem \ref{normquantile1}.}
Since
\begin{eqnarray*}
F(\mu_p(x)|x)=p=F_\nbf(\mu_{p,\nbf}(x)|x),
\end{eqnarray*}
by Taylor's expansion, we have
\begin{eqnarray*}
\mu_{p,\nbf}(x)-\mu_p(x)=\frac{1}{f_\nbf(\mu_{p,\nbf}^{*}(x)|x)}\left[F(\mu_p(x)|x)-F_{\nbf}(\mu_p(x)|x)\right],
\end{eqnarray*}
where $\mu_{p,\nbf}^{*}(x)$ lies between $\mu_{p,\nbf}(x)$ and $\mu_p(x)$. To prove the asymptotic normality, it suffices to show that the numerator is normally distributed, and the denominator converges to $f(\mu_p(x)|x)$ in probability. We have the following propositions.
\begin{prop}\label{normcond}
Under (A1)-(A6), (C1) and (C2), if there exists $c\geq 0$ such that $\widehat\nbf h^{d+4}\to c$, then
\begin{eqnarray*}
\sqrt{\widehat\nbf h^d}\left(F_{\nbf}(y|x)-F(y|x)\right)\stackrel{\mathcal{L}}{\to} \mathcal{N}\left(cB(x,y),\sigma^2(x,y)\right),
\end{eqnarray*}
where $B(x,y)$ and $\sigma^2(x,y)$ are defined in \eqref{biais} and \eqref{variance}.
\end{prop}
\begin{prop}\label{convunifconddens}
Assume that (A1), (A3)-(A7) hold. If (A9) and (\ref{mixing1}) or (A10) and (\ref{mixing2}) are satisfied, then
\begin{eqnarray*}
\sqrt{\frac{\widehat\nbf h^{d+1}}{\log \widehat\nbf}}\sup_{y\in \mathcal{V}}|f_\nbf(y|x)-f(y|x)| \stackrel{a.s.}{\to} 0.
\end{eqnarray*}
\end{prop}
\noindent\textbf{Proof of Proposition \ref{normcond}.}
Proposition \ref{normcond} is a consequence of the two following lemmas.
\begin{lem}\label{condbias}
If Assumptions (A1) and (A3)-(A6) are satisfied, then
\begin{eqnarray*}
\Ebf F_{\nbf}(y|x)-F(y|x)=\frac{h^2}{2}B(x,y)+o(h^2)+O\left(\frac{1}{\widehat\nbf h^d}\right),
\end{eqnarray*}
where $B(x,y)$ is defined in \eqref{biais}.
\end{lem}
The proof of Lemma \ref{condbias} is classical and therefore is omitted (see Matzner-L\o ber \cite{matzner}).
\begin{lem}\label{normcondvar}
Assume that Assumptions (A1)-(A6), (C1) and (C2) hold, then
\begin{eqnarray*}
\sqrt{\widehat\nbf h^d}\left(F_{\nbf}(y|x)-\Ebf F_{\nbf}(y|x)\right)\stackrel{\mathcal{L}}{\to} \mathcal{N}\left(0,\sigma^{2}(x,y)\right),
\end{eqnarray*}
 where $\sigma^2(x,y)$ is defined in \eqref{variance}.
\end{lem}
\noindent\textbf{Proof of Lemma \ref{normcondvar}.}
Assume for the moment that
 for any pair $(c_1,c_2) \in \mathbb{R}^2$ with $c_1^2+c_2^2\neq 0$,
\begin{eqnarray}
\label{norm}
\sqrt{\widehat\nbf h^d}\left[c_1\left(g_\nbf(x)-\Ebf g_\nbf(x)\right)+c_2\left(\psi_{\nbf}(x,y)-\Ebf\psi_{\nbf}(x,y)\right)\right]\stackrel{\mathcal{L}}{\to}\mathcal{N}\left(0,\sigma^*{^2}\right),
\end{eqnarray}
where $\displaystyle \sigma^*{^2}= \left[c_1^2g(x)+c_2^2\psi(x,y)+2c_1c_2\psi(x,y)\right]\int_{\mathbb{R}^d}K^2(z)dz.$ Now, set
\begin{eqnarray*}
W_\nbf &=&\frac{1}{g_\nbf(x)}\left[\psi_{\nbf}(x,y)-\Ebf\psi_{\nbf}(x,y)\right]-\frac{\Ebf\psi_{\nbf}(x,y)}{g_\nbf(x)\Ebf g_\nbf(x)}\left[g_\nbf(x)-\Ebf g_\nbf(x)\right],\\
W_\nbf^*&=&\frac{1}{g(x)}\left[\psi_{\nbf}(x,y)-\Ebf\psi_{\nbf}(x,y)\right]-\frac{\psi(x,y)}{g^2(x)}\left[g_\nbf(x)-\Ebf g_\nbf(x)\right].
\end{eqnarray*}
On one hand, according to \eqref{norm}, we have
$\sqrt{\widehat\nbf h^d}W_\nbf^*\stackrel{\mathcal{L}}{\to}\mathcal{N}\left(0,\sigma{^2}\right)$,
with $\sigma^2$ being defined in \eqref{variance}.\\
On the other hand,
\begin{eqnarray*}
\lefteqn{\sqrt{\widehat\nbf h^d}\left[W_\nbf-W_\nbf^*\right]}\\
 &= & \sqrt{\widehat\nbf h^d}\left\{\left[\frac{1}{g_\nbf(x)}-\frac{1}{g(x)}\right]\left[\psi_{\nbf}(x,y)-\Ebf \psi_{\nbf}(x,y)\right]+\left[\frac{\psi(x,y)}{g^2(x)}-\frac{\Ebf\psi_{\nbf}(x,y)}{g_\nbf(x)\Ebf g_\nbf(x)}\right]\left[g_\nbf(x)-\Ebf g_\nbf(x)\right]\right\}
\end{eqnarray*}
which goes to $0$ in probability. \\
Finally, we conclude with the decomposition
\begin{eqnarray*}
\sqrt{\widehat\nbf h^d}\left(F_{\nbf}(y|x)-\Ebf F_{\nbf}(y|x)\right)= \sqrt{\widehat\nbf h^d}W_\nbf^*+\sqrt{\widehat\nbf h^d}\left[W_\nbf-W_\nbf^*\right]+O\left(\frac{1}{\sqrt{\widehat\nbf h^d}}\right).
\end{eqnarray*}
Now, \eqref{norm} is proved following the same lines as the proof of Theorem 3.1 in Tran \cite{tran}.\\

\noindent\textbf{Proof of Proposition \ref{convunifconddens}.}
We have
\begin{eqnarray*}
\sup_{y \in \mathcal{V}}|f_\nbf(y|x)-f(y|x)|&\leq & \frac{1}{g_\nbf(x)}\sup_{y \in \mathcal{V}}|f_\nbf(x,y)-f(x,y)|+\frac{1}{g_\nbf(x)}\sup_{y \in \mathcal{V}}f(y|x)|g_\nbf(x)-g(x)|.
\end{eqnarray*}
Hence to prove Proposition \ref{convunifconddens}, it suffices to show that
$\displaystyle\sup_{y \in \mathcal{V}}|f_\nbf(x,y)-f(x,y)|\stackrel{\Pbf}{\to}0$, this is proved by following the same lines as in the proof of Theorem 3.3 in Carbon et al. \cite{carbon}.

\end{document}